\documentclass[10pt,portrait,reqno,12pt]{amsart}%
\usepackage{amssymb,amsthm,amsmath, amsfonts, float}
\usepackage{graphicx}
\usepackage{multirow}
\usepackage{subfig}
\usepackage{amsmath}
\usepackage{amsfonts}
\usepackage{lscape}
\usepackage{amssymb}%
\providecommand{\U}[1]{\protect\rule{.1in}{.1in}}
\newtheorem{theorem}{Theorem}
\theoremstyle{plain}

\newtheorem{definition}{Definition}

\newtheorem{lemma}{Lemma}

\numberwithin{equation}{section}
\begin{document}
\title[$\mathcal{L}_{1}$ AND $\mathcal{L}_{2}$ OPERATORS\ OF\ GAUSS MAP FOR TUBULAR
HYPERSURFACES IN $\mathbb{E}%
^{4}_{1}$]{SOME\ CLASSIFICATIONS FOR GAUSS MAP OF TUBULAR HYPERSURFACES IN $\mathbb{E}%
^{4}_{1}$
CONCERNING LINEARIZED OPERATORS\ $\mathcal{L}_{k}$}
\author[A. Kazan, M. Alt\i n and N.C. Turgay]{\bfseries Ahmet Kazan$^{1\ast}$, Mustafa Alt\i n$^{2}$ and Nurettin Cenk
Turgay$^{3}$}
\address{
$^{1}$\textit{Department of Computer Technologies, Do\u{g}an\c{s}ehir Vahap
K\"{u}\c{c}\"{u}k Vocational School, Malatya Turgut \"{O}zal University,
Malatya, Turkey}
\newline
$^{2}$\textit{Technical Sciences Vocational School, Bing\"{o}l University,
Bing\"{o}l, Turkey}
\newline
$^{3}$\textit{Department of Mathematics, Faculty of Science and Letters,
\.{I}stanbul Technical University, \.{I}stanbul, Turkey}
\newline
$^{\ast}$\textit{Corresponding author: ahmet.kazan@ozal.edu.tr}}

\begin{abstract}
In this study, we deal with the Gauss map of tubular hypersurfaces in
4-dimensional Lorentz-Minkowski space concerning the linearized operators
$\mathcal{L}_{1}$ (Cheng-Yau) and $\mathcal{L}_{2}$. We obtain the
$\mathcal{L}_{1}$ (Cheng-Yau) operator of the Gauss map of tubular
hypersurfaces that are formed as the envelope of a family of pseudo
hyperspheres{ or pseudo hyperbolic hyperspheres} whose centers lie on timelike
or spacelike curves with non-null Frenet vectors in $\mathbb{E}%
^{4}_{1}$ and give some
classifications for these hypersurfaces which have generalized $\mathcal{L}%
_{k}$ 1-type Gauss map, first kind $\mathcal{L}_{k}$-pointwise 1-type Gauss
map, second kind $\mathcal{L}_{k}$-pointwise 1-type Gauss map and $\mathcal{L}_{k}$-harmonic Gauss map, $k\in\{1,2\}$.

\end{abstract}
\maketitle

%=====Contents=======

\section{Introduction}

\label{SectionIntrod} Let $(M,g)$ be a hypersurface of $(n+1)$-dimensional
Minkowski space $\mathbb{E}^{n+1}_{1}$, $\Delta$ denote its Laplace operator.
A smooth mapping $\phi:M\rightarrow\mathbb{E}^{n+1}_{1}$ is said to be \textit{finite type} if it can be expressed as
\[
\phi=\phi_{0}+\phi_{1}+\cdots+\phi_{k},
\]
where $\phi_{0}$ is a constant vector and $\phi_{i}$ is an eigenvector of
$\Delta$ corresponding to the eigenvector $\lambda_{i}$ for $i=1,2,\hdots,k$.
More precisely, if $\lambda_{1},\lambda_{2},\hdots,\lambda_{k}$ are distinct,
then $\psi$ is said to be $k$-\textit{type}
(\cite{ChenKitap,Chen-Morvan-Nore,Chen-Petru}). Several results on the study
of finite type mappings were summed up in a report by B.-Y. Chen in
\cite{ChenRapor} (See also \cite{Chen-Piccinni,YoonTorus1}).

Let $N$ denote the Gauss map of $M$. From the definition above, one can
conclude that $N$ is of \textit{1-type} if and only if it satisfies the equation
\begin{equation}
\label{Glob1TpeDef}\Delta N = \lambda(N+C)
\end{equation}
for a constant $\lambda\in\mathbb{R}$ and a constant vector $C$. However,
Gauss map of some important submanifolds such as catenoid and helicoid of the
Euclidean 3-space $\mathbb{E}^{3}$ satisfies
\begin{equation}
\label{PW1TypeDefinition}\Delta N =f(N+C)
\end{equation}
which is weaker than \eqref{Glob1TpeDef}, where $f\in C^{\infty}(M)$
is a smooth function, \cite{ChoiKim2001}. These submanifolds whose
Gauss map $N$ satisfying \eqref{PW1TypeDefinition} are said to have
\textit{pointwise 1-type Gauss map}. Submanifolds with pointwise
1-type Gauss map have been worked in several papers (cf.
\cite{ChoiKim2001,SenYegin2017,Takahashi,YoonTorus1}).

On the other hand, the Gauss map of some hypersurfaces of semi-Euclidean
spaces satisfies the equation
\begin{equation}
\Delta N=f_{1}N+f_{2}C \label{03}%
\end{equation}
for some smooth functions $f_{1}$, $f_{2}$ and a constant vector $C$. A
submanifold is said to have \textit{generalized 1-type Gauss map} if its Gauss
map satisfies the condition \eqref{03}, \cite{yoon2018}. After this definition
was given, hypersurfaces of pseudo-Euclidean spaces have been considered in
terms of having generalized 1-type Gauss map,
\cite{qian2020,kim2021,yoon2018,yoon2018b}.

In the recent years, the definition of $\mathcal{L}_{k}$-finite type maps has
been obtained by replacing $\Delta$ in the definition above with the sequence
of operators $\mathcal{L}_{0},\ \mathcal{L}_{1}, \ \mathcal{L}_{2}%
,\ \hdots,\ \mathcal{L}_{n-1}$, \cite{Alias2,AliasGurbuzGeomDed2006}. Note
that, by the definition of these operators, one can obtain $\mathcal{L}%
_{0}=-\Delta$ and $\mathcal{L}_{1}=\square$ is called as the
Cheng-Yau operator introduced in \cite{Cheng-Yau}. By motivating
this idea, notion of $\mathcal{L}_{k}$-pointwise 1-type Gauss map
and generalized $\mathcal{L}_{k}$ 1-type Gauss map was presented in
\cite{KimTurgay} and \cite{qian2020b}, respectively (see Definition
\ref{MainDef}). After the case $k=1$ is studied in these papers,
many result obtained on hypersurfaces with certain type of Gauss
map, \cite{Kazany,Kelleci,Kim2,kim2021,Yang2,yoon2018,yoon2018b}.

On the other hand, in \cite{BKMSbizim}, the general expression of the canal
hypersurfaces that are formed as the envelope of a family of pseudo
hyperspheres, pseudo hyperbolic hyperspheres and null hypercones whose centers
lie on a non-null curve with non-null Frenet vector fields in $\mathbb{E}%
^{4}_{1}$ has
been given and their some geometric invariants such as unit normal vector
felds, Gaussian curvatures, mean curvatures and principal curvatures have been
obtained. Also, tubular hypersurfaces in $\mathbb{E}%
^{4}_{1}$ by taking constant radius function have been studied in
\cite{BKMSbizim}.

In this paper, we study the tubular hypersurfaces in Lorentz-Minkowski 4-space $\mathbb{E}%
^{4}_{1}$ with the aid of $\mathcal{L}_{k}$ operators, $k\in\{1,2\}$. In Sect. 2, we give basic
notation, facts and definitions about hypersurfaces of Minkowski spaces. In
Sect. 3 and Sect 4, we consider some classifications of tubular hypersurfaces by
considering their Gauss maps in terms of their types with respect to the
operators $\mathcal{L}_{1}$ and $\mathcal{L}_{2}$.

\section{{Preliminaries}}

Let $\mathbb{\mathbb{E}}_{1}^{n+1}$ be the $(n+1)$-dimensional
Lorentz-Minkowski space with the canonical pseudo-Euclidean metric
$\left\langle \text{ },\right\rangle $ of index 1 and signature
$(-,+,+,...,+)$ given by%
\[
\left\langle \text{ },\right\rangle =-dx_{1}^{2}+dx_{2}^{2}+dx_{3}%
^{2}+...+dx_{n+1}^{2}%
\]
where $(x_{1},x_{2},...,x_{n+1})$ is a rectangular coordinate system in
$\mathbb{\mathbb{E}}_{1}^{n+1}$.

If $\Gamma:M\longrightarrow\mathbb{E}_{1}^{n+1}$ is an isometric immersion
from an $n$-dimensional orientable manifold $M$ to $\mathbb{\mathbb{E}}%
_{1}^{n+1}$, then the induced metric on $M$ by the immersion $\Gamma$ can be
Riemannian or Lorentzian. Let $N$ denotes a unit normal vector field and put
$\left\langle N,N\right\rangle =\varepsilon=\pm1$, so that $\varepsilon=1$ or
$\varepsilon=-1$ according to $M$ is endowed with a Lorentzian or Riemannian
metric, respectively.

The operator $\mathcal{L}_{k}$ acting on the coordinate functions of the Gauss
map $N$ of the hypersurface $M$ in $(n+1)$-dimensional Lorentz-Minkowski space
$\mathbb{\mathbb{E}}_{1}^{n+1}$ is%
\begin{equation}
\mathcal{L}_{k}N=-\varepsilon\mathfrak{C}_{k}\left(  \nabla H_{k+1}+\left(
nH_{1}H_{k+1}-\left(  n-k-1\right)  H_{k+2}\right)  N\right)  . \label{Lk}%
\end{equation}
Here,%
\begin{equation}
\left(
\begin{array}
[c]{c}%
n\\
k
\end{array}
\right)  H_{k}=\left(  -\varepsilon\right)  ^{k}a_{k},\text{ \ }\left(
\left(
\begin{array}
[c]{c}%
n\\
k
\end{array}
\right)  =\frac{n!}{k!(n-k)!}\right)  , \label{Hk}%
\end{equation}
such that%

\begin{equation}
\left.
\begin{array}
[c]{l}%
a_{1}=-\sum_{i=1}^{n}\kappa_{i},\\
\\
a_{k}=(-1)^{k}\sum_{i_{1}<i_{2}<...<i_{k}}^{n}\kappa_{i_{1}}\kappa_{i_{2}%
}...\kappa_{i_{k}},~k=2,3,...,n
\end{array}
\right\}  \label{ak}%
\end{equation}
and $H_{k}$ is called the $k$-th mean curvature of order $k$ of $M$.

Also, the constant $\mathfrak{C}_{k}$ is given by%
\begin{equation}
\mathfrak{C}_{k}=\left(
\begin{array}
[c]{c}%
n\\
k+1
\end{array}
\right)  (-\varepsilon)^{k}. \label{Ck}%
\end{equation}
(For more details about the linearized operator $\mathcal{L}_{k}$, one can see
\cite{ChengYauEn1}.)

\begin{definition}
\label{def}Let $\mathfrak{m}$ and $\mathfrak{n}$ be non-zero smooth functions
on $M$, $C\in\mathbb{E}_{1}^{n+1}$ be a non-zero constant vector and
$k\in\{0,1,2,...,n\}$.

If the Gauss map $N$ of an oriented submanifold $M$ in $\mathbb{E}%
^{4}_{1}$ satisfies

\begin{description}
\label{MainDef}

\item[i] $\mathcal{L}_{k}N=\mathfrak{m}N+\mathfrak{n}C,$ then $M$ has
generalized $\mathcal{L}_{k}$ 1-type Gauss map;

\item[ii] $\mathcal{L}_{k}N=\mathfrak{m}N,$ then $M$ has first kind
$\mathcal{L}_{k}$-pointwise 1-type Gauss map;

\item[iii] $\mathcal{L}_{k}N=\mathfrak{m}(N+C),$ then $M$ has second kind
$\mathcal{L}_{k}$-pointwise 1-type Gauss map;

\item[iv] $\mathcal{L}_{k}N=0$, then $N$ is called $\mathcal{L}_{k}$-harmonic.
\end{description}
\end{definition}

In this study, we will deal with Gauss maps of tubular hypersurfaces in
4-dimensional Lorentz-Minkowski space $\mathbb{E}%
^{4}_{1}$ concerning linearized
operators $L_{1}$ and $L_{2}$. So, let us give some notions in $\mathbb{E}%
^{4}_{1}$.

Let $\overrightarrow{u}=(u_{1},u_{2},u_{3},u_{4})$, $\overrightarrow{v}%
=(v_{1},v_{2},v_{3},v_{4})$ and $\overrightarrow{w}=(w_{1},w_{2},w_{3},w_{4})$
be three vectors in $\mathbb{E}%
^{4}_{1}$. The inner product and vector product are
defined by%
\begin{equation}
\left\langle \overrightarrow{u},\overrightarrow{v}\right\rangle =-u_{1}%
v_{1}+u_{2}v_{2}+u_{3}v_{3}+u_{4}v_{4} \label{yy1}%
\end{equation}
and
\begin{equation}
\overrightarrow{u}\times\overrightarrow{v}\times\overrightarrow{w}=\det\left[
\begin{array}
[c]{cccc}%
-e_{1} & e_{2} & e_{3} & e_{4}\\
u_{1} & u_{2} & u_{3} & u_{4}\\
v_{1} & v_{2} & v_{3} & v_{4}\\
w_{1} & w_{2} & w_{3} & w_{4}%
\end{array}
\right]  , \label{yy2}%
\end{equation}
respectively. Here $e_{i},(i=1,2,3,4)$ are standard basis vectors.

A vector $\overrightarrow{u}\in E_{1}^{4}-\{0\}$ is called spacelike, timelike
or lightlike (null) if $\left\langle \overrightarrow{u},\overrightarrow{u}%
\right\rangle >0$ (or $\overrightarrow{u}=0$), $\left\langle
\overrightarrow{u},\overrightarrow{u}\right\rangle <0$ or $\left\langle
\overrightarrow{u},\overrightarrow{u}\right\rangle =0$, respectively. A curve
$\beta(s)$ in $\mathbb{E}%
^{4}_{1}$ is spacelike, timelike or lightlike (null), if all
its velocity vectors $\beta^{\prime}(s)$ are spacelike, timelike or lightlike,
respectively and a non-null (i.e. timelike or spacelike) curve has unit speed
if $\left\langle \beta^{\prime},\beta^{\prime}\right\rangle =\mp1$. Also, the
norm of the vector $\overrightarrow{u}$ is $\left\Vert \overrightarrow{u}%
\right\Vert =\sqrt{\left\vert \left\langle \overrightarrow{u}%
,\overrightarrow{u}\right\rangle \right\vert }$ \cite{Kuhnel}.

Let $F_{1},$ $F_{2},$ $F_{3},$ $F_{4}$ be unit tangent vector field, principal
normal vector field, binormal vector field, trinormal vector field of a
timelike or spacelike curve $\beta(s)$, respectively and $\{F_{1},F_{2}%
,F_{3},F_{4}\}$ be the moving Frenet frame along $\beta(s)$ in $\mathbb{E}%
^{4}_{1}$.
The Frenet equations can be given according to the causal characters of
non-null Frenet vector fields $F_{1},$ $F_{2},$ $F_{3}$ and $F_{4}$ as follows
\cite{Walrave}:

If the curve $\beta(s)$ is timelike, i.e. $\left\langle F_{1},F_{1}%
\right\rangle =-1,$ $\left\langle F_{i},F_{i}\right\rangle =1$ ($i=2,3,4$)$,$
then%
\begin{equation}
\left.
\begin{array}
[c]{l}%
F_{1}^{\prime}=k_{1}F_{2},\\
F_{2}^{\prime}=k_{1}F_{1}+k_{2}F_{3},\\
F_{3}^{\prime}=-k_{2}F_{2}+k_{3}F_{4},\\
F_{4}^{\prime}=-k_{3}F_{3};
\end{array}
\right\}  \label{fr1}%
\end{equation}
if the curve $\beta(s)$ is spacelike with timelike principal normal vector
field $F_{2}$, i.e. $\left\langle F_{2},F_{2}\right\rangle =-1,$ $\left\langle
F_{i},F_{i}\right\rangle =1$ ($i=1,3,4$)$,$ then%
\begin{equation}
\left.
\begin{array}
[c]{l}%
F_{1}^{\prime}=k_{1}F_{2},\\
F_{2}^{\prime}=k_{1}F_{1}+k_{2}F_{3},\\
F_{3}^{\prime}=k_{2}F_{2}+k_{3}F_{4},\\
F_{4}^{\prime}=-k_{3}F_{3};
\end{array}
\right\}  \label{fr2}%
\end{equation}
if the curve $\beta(s)$ is spacelike with timelike binormal vector field
$F_{3}$, i.e. $\left\langle F_{3},F_{3}\right\rangle =-1,$ $\left\langle
F_{i},F_{i}\right\rangle =1$ ($i=1,2,4$)$,$ then%
\begin{equation}
\left.
\begin{array}
[c]{l}%
F_{1}^{\prime}=k_{1}F_{2},\\
F_{2}^{\prime}=-k_{1}F_{1}+k_{2}F_{3},\\
F_{3}^{\prime}=k_{2}F_{2}+k_{3}F_{4},\\
F_{4}^{\prime}=k_{3}F_{3};
\end{array}
\right\}  \label{fr3}%
\end{equation}
if the curve $\beta(s)$ is spacelike with timelike trinormal vector field
$F_{4}$, i.e. $\left\langle F_{4},F_{4}\right\rangle =-1,$ $\left\langle
F_{i},F_{i}\right\rangle =1$ ($i=1,2,3$)$,$ then%
\begin{equation}
\left.
\begin{array}
[c]{l}%
F_{1}^{\prime}=k_{1}F_{2},\\
F_{2}^{\prime}=-k_{1}F_{1}+k_{2}F_{3},\\
F_{3}^{\prime}=-k_{2}F_{2}+k_{3}F_{4},\\
F_{4}^{\prime}=k_{3}F_{3}.
\end{array}
\right\}  \label{fr4}%
\end{equation}
Here $k_{1},k_{2},k_{3}$ are the first, second and third curvatures of the
non-null curve $\beta(s).$

Also, if $p$ is a fixed point in $\mathbb{E}%
^{4}_{1}$ and $r$ is a positive constant,
then the pseudo-Riemannian hypersphere and the pseudo-Riemannian hyperbolic
space are defined by%
\[
S_{1}^{3}(p,r)=\{x\in \mathbb{E}%
^{4}_{1}:\left\langle x-p,x-p\right\rangle =r^{2}\}
\]
and
\[
H_{0}^{3}(p,r)=\{x\in \mathbb{E}%
^{4}_{1}:\left\langle x-p,x-p\right\rangle =-r^{2}\},
\]
respectively.

If $M$ is an oriented hypersurface in $E_{1}^{4},$ then the gradient of a
smooth function $f(s,t,w)$, which is defined in $M$, can be obtained by%
\begin{equation}
\nabla f=\frac{1}{\mathfrak{g}}\left(
\begin{array}
[c]{c}%
\left(  \left(  g_{23}^{2}-g_{22}g_{33}\right)  f_{s}+\left(  -g_{13}%
g_{23}+g_{12}g_{33}\right)  f_{t}+\left(  g_{13}g_{22}-g_{12}g_{23}\right)
f_{w}\right)  \partial s\\
+\left(  \left(  -g_{13}g_{23}+g_{12}g_{33}\right)  f_{s}+\left(  g_{13}%
^{2}-g_{11}g_{33}\right)  f_{t}+\left(  -g_{12}g_{13}+g_{11}g_{23}\right)
f_{w}\right)  \partial t\\
+\left(  \left(  g_{13}g_{22}-g_{12}g_{23}\right)  f_{s}+\left(  -g_{12}%
g_{13}+g_{11}g_{23}\right)  f_{t}+\left(  g_{12}^{2}-g_{11}g_{22}\right)
f_{w}\right)  \partial w
\end{array}
\right)  , \label{gradf}%
\end{equation}
where%
\[
\mathfrak{g=}g_{13}^{2}g_{22}-2g_{12}g_{13}g_{23}+g_{11}g_{23}^{2}+g_{12}%
^{2}g_{33}-g_{11}g_{22}g_{33};
\]
$\{s,t,w\}$ is a local coordinat system of $M$; $f_{s}$, $f_{t}$, $f_{w}$ are
the partial derivatives of $f$ and $g_{11}=\left\langle \partial s,\partial
s\right\rangle ,$ $g_{12}=\left\langle \partial s,\partial t\right\rangle ,$
$g_{13}=\left\langle \partial s,\partial w\right\rangle ,$ $g_{22}%
=\left\langle \partial t,\partial t\right\rangle ,$ $g_{23}=\left\langle
\partial t,\partial w\right\rangle ,$ $g_{33}=\left\langle \partial w,\partial
w\right\rangle $.

\section{Some Classifications for Tubular Hypersurfaces Generated by Timelike Curves with $L_{k}$ Operators in $\mathbb{E}%
^{4}_{1}$}

In this section, we obtain the $\mathcal{L}_{1}$ (Cheng-Yau) and
$\mathcal{L}_{2}$ operators of the Gauss map of the tubular hypersurfaces
$\mathcal{T}(s,t,w)$ that are formed as the envelope of a family of pseudo
hyperspheres whose centers lie on a timelike curve with non-null Frenet
vectors in $\mathbb{E}%
^{4}_{1}$ and give some classifications for these hypersurfaces
which have generalized $\mathcal{L}_{k}$ 1-type Gauss map,\ first kind
$\mathcal{L}_{k}$-pointwise 1-type Gauss map and second kind $\mathcal{L}_{k}%
$-pointwise 1-type Gauss map and $\mathcal{L}_{k}$-harmonic Gauss map,
$k\in\{1,2\}$.

The tubular hypersurfaces $\mathcal{T}(s,t,w)$ that are formed as the envelope
of a family of pseudo hyperspheres whose centers lie on a timelike curve with
non-null Frenet vectors in $\mathbb{E}%
^{4}_{1}$ can be parametrized by%
\begin{equation}
\mathcal{T}(s,t,w)=\beta(s)+r(\cos t\cos wF_{2}(s)+\sin t\cos wF_{3}(s)+\sin
wF_{4}(s)). \label{surfT11}%
\end{equation}
The unit normal vector field of {(\ref{surfT11}) is}%
\begin{equation}
N={-}\left(  {\cos t\cos wF}_{2}+{\sin t\cos wF}_{3}+{\sin wF}_{4}\right)
\label{normalT11}%
\end{equation}
and so,
\begin{equation}
\left\langle N,N\right\rangle =1. \label{epsilonT11}%
\end{equation}
The coefficients of the first fundamental form of {(\ref{surfT11}) are}%
\begin{equation}
\left.
\begin{array}
[c]{l}%
g_{11}={\small -(1+rk}_{1}\cos{\small t}\cos{\small w)}^{2}{\small +(rk}%
_{2}\cos{\small t}\cos{\small w-rk}_{3}\sin{\small w)}^{2}{\small +r}%
^{2}\left(  {\small k}_{2}^{2}{\small +k}_{3}^{2}\right)  \sin^{2}%
{\small t}\cos^{2}{\small w,}\\
g_{12}=g_{21}=r^{2}(k_{2}\cos w-k_{3}\cos t\sin w)\cos w,\text{ }g_{22}%
=r^{2}\cos^{2}w,\\
g_{13}=g_{31}=r^{2}k_{3}\sin t,\text{ }g_{23}=g_{32}=0,\text{ }g_{33}=r^{2}.
\end{array}
\right\}  \label{gT11}%
\end{equation}
The principal curvatures of {(\ref{surfT11}) are}%
\begin{equation}
\kappa_{1}=\kappa_{2}=\frac{1}{r},\text{ }\kappa_{3}=\frac{k_{1}\cos t\cos
w}{1+rk_{1}\cos t\cos w}. \label{asliT11}%
\end{equation}
For more details about these hypersurfaces, one can see \cite{BKMSbizim}.

\subsection{Some Classifications for Tubular Hypersurfaces Generated by Timelike Curves with $\mathcal{L}_{1}$ (Cheng-Yau) Operator in $\mathbb{E}%
    ^{4}_{1}
$}

\

The functions $a_{k}$ of the tubular hypersurfaces (\ref{surfT11}) in
$\mathbb{E}%
^{4}_{1}$ are obtained from (\ref{ak}) and (\ref{asliT11}) by%
\begin{equation}
a_{1}=\frac{-2-3rk_{1}\cos t\cos w}{r(1+rk_{1}\cos t\cos w)},\text{ }%
a_{2}=\frac{1+3rk_{1}\cos t\cos w}{r^{2}(1+rk_{1}\cos t\cos w)},\text{ }%
a_{3}=-\frac{k_{1}}{r^{2}(rk_{1}+\sec t\sec w)}. \label{akT11}%
\end{equation}
Also, from (\ref{gradf}), (\ref{gT11}) and (\ref{akT11}), we have%
\begin{align}
\nabla a_{2}  &  ={-\frac{2\left(  k_{1}^{\prime}\cos t+k_{1}k_{2}\sin
t\right)  \cos w}{r(1+rk_{1}\cos t\cos w)^{3}}}F_{1}{-\frac{k_{1}\left(
2\cos^{2}t\cos(2w)+\cos(2t)-3\right)  }{2r^{2}(1+rk_{1}\cos t\cos w)^{2}}%
}F_{2}\nonumber\\
&  {-\frac{2k_{1}\sin t\cos t\cos^{2}w}{r^{2}(1+rk_{1}\cos t\cos w)^{2}}}%
F_{3}{-\frac{2k_{1}\cos t\sin w\cos w}{r^{2}(1+rk_{1}\cos t\cos w)^{2}}}F_{4}.
\label{a2nablaT11}%
\end{align}

So, from (\ref{Lk}), (\ref{Hk}), (\ref{normalT11}), (\ref{epsilonT11}),
(\ref{akT11}) and (\ref{a2nablaT11}), we reach that%
\begin{align}
\mathcal{L}_{1}N  &  ={-\frac{2\left(  k_{1}k_{2}\sin t+k_{1}^{\prime}\cos
t\right)  \cos w}{r(1+rk_{1}\cos t\cos w)^{3}}}F_{1}\nonumber\\
&  {-\frac{2\left(  rk_{1}\left(  3rk_{1}\cos^{3}t\cos^{3}w+2\cos^{2}%
t\cos(2w)+\cos(2t)\right)  +\cos t\cos w\right)  }{r^{3}(1+rk_{1}\cos t\cos
w)^{2}}}F_{2}\nonumber\\
&  {-\frac{2(1+3rk_{1}\cos t\cos w)\sin t\cos w}{r^{3}(1+rk_{1}\cos t\cos w)}%
}F_{3}{-\frac{2(3rk_{1}\cos t\cos w+1)\sin w}{r^{3}(1+rk_{1}\cos t\cos w)}%
}F_{4}. \label{L1T11}%
\end{align}

Now, let us give some classifications for the tubular hypersurfaces
(\ref{surfT11}) which have generalized $\mathcal{L}_{1}$ 1-type Gauss
map,\ first kind $\mathcal{L}_{1}$-pointwise 1-type Gauss map and second kind
$\mathcal{L}_{1}$-pointwise 1-type Gauss map and $\mathcal{L}_{1}$-harmonic Gauss map.

Let the tubular hypersurfaces $\mathcal{T}(s,t,w)$ have generalized
$\mathcal{L}_{1}$ (Cheng-Yau) 1-type Gauss map, i.e., $\mathcal{L}%
_{1}N=\mathfrak{m}N+\mathfrak{n}C,$ where $C=C_{1}F_{1}+C_{2}F_{2}+C_{3}%
F_{3}+C_{4}F_{4}$ is a constant vector. Here, by taking derivatives of the
constant vector $C$ with respect to $s,$ from (\ref{fr1}) we obtain that%
\begin{equation}
\left.
\begin{array}
[c]{l}%
C_{1}^{\prime}+C_{2}k_{1}=0,\\
C_{2}^{\prime}+C_{1}k_{1}-C_{3}k_{2}=0,\\
C_{3}^{\prime}+C_{2}k_{2}-C_{4}k_{3}=0,\\
C_{4}^{\prime}+C_{3}k_{3}=0.
\end{array}
\right\}  \label{turevT11}%
\end{equation}
Also, by taking derivatives the constant vector $C$ with respect to $t$ and $w$
separately, one can see that the functions $C_{i}$ depend\ only on $s.$

Firstly, let us classificate the tubular hypersurfaces $\mathcal{T}(s,t,w)$
which have generalized $\mathcal{L}_{1}$ (Cheng-Yau) 1-type Gauss map.

From (\ref{normalT11}) and (\ref{L1T11}), we get%
\begin{equation}
\left.
\begin{array}
[c]{l}%
{-\frac{2\left(  k_{1}^{\prime}\cos t+k_{1}k_{2}\sin t\right)  \cos
w}{r(1+rk_{1}\cos t\cos w)^{3}}}=\mathfrak{n}C_{1},\\
\\
{-\frac{2\left(  rk_{1}\left(  3rk_{1}\cos^{3}t\cos^{3}w+2\cos^{2}%
t\cos(2w)+\cos(2t)\right)  +\cos t\cos w\right)  }{r^{3}(1+rk_{1}\cos t\cos
w)^{2}}}=\mathfrak{m}\left(  {-\cos t\cos w}\right)  +\mathfrak{n}C_{2},\\
\\
{-\frac{2(1+3rk_{1}\cos t\cos w)\sin t\cos w}{r^{3}(1+rk_{1}\cos t\cos w)}%
}=\mathfrak{m}\left(  -{\sin t\cos w}\right)  +\mathfrak{n}C_{3},\\
\\
{-\frac{2(1+3rk_{1}\cos t\cos w)\sin w}{r^{3}(1+rk_{1}\cos t\cos w)}%
}=\mathfrak{m}\left(  -{\sin w}\right)  +\mathfrak{n}C_{4}.
\end{array}
\right\}  \label{denklemlerT11}%
\end{equation}

Now let us investigate the non-zero functions $\mathfrak{m}(s,t,w)$ and
$\mathfrak{n}(s,t,w)$ from the above four equations.

Firstly, let us assume that $C_{1}\neq0$.

In this case, from the first equation of (\ref{denklemlerT11}) it's easy to
see that%
\begin{equation}
\mathfrak{n}(s,t,w)={-\frac{2\left(  k_{1}^{\prime}\cos t+k_{1}k_{2}\sin
t\right)  \cos w}{r(1+rk_{1}\cos t\cos w)^{3}C_{1}}.} \label{c1esitdegil0nT11}%
\end{equation}
Here, when the equation (\ref{c1esitdegil0nT11}) is successively substituted
into the second, third and fourth equations of (\ref{denklemlerT11}), we
obtain%
\[
\left.
\begin{array}
[c]{l}%
\mathfrak{m}(s,t,w)=\frac{2\left(
\begin{array}
[c]{l}%
\left(  {\small rk}_{1}\left(
\begin{array}
[c]{l}%
{\small 3rk}_{1}\cos^{3}{\small t}\cos^{3}{\small w}\\
{\small +2}\cos^{2}{\small t}\cos{\small (2w)+}\cos{\small (2t)}%
\end{array}
\right)  {\small +}\cos{\small t}\cos{\small w}\right)  {\small C}%
_{1}{\small (1+rk}_{1}\cos{\small t}\cos{\small w)}\\
\\
{\small -C}_{2}{\small r}^{2}\left(  k_{1}^{\prime}\cos t+k_{1}k_{2}\sin
t\right)  \cos w
\end{array}
\right)  }{C_{1}r^{3}(1+rk_{1}\cos t\cos w)^{3}\cos t\cos w},\\
\\
\mathfrak{m}(s,t,w)=\frac{2\left(  C_{1}(1+rk_{1}\cos t\cos w)^{2}%
(1+3rk_{1}\cos t\cos w)-C_{3}r^{2}\left(  k_{1}^{\prime}\cot t+k_{1}%
k_{2}\right)  \right)  }{C_{1}r^{3}(1+rk_{1}\cos t\cos w)^{3}},\\
\\
\mathfrak{m}(s,t,w)=\frac{2\left(  C_{1}(1+rk_{1}\cos t\cos w)^{2}%
(1+3rk_{1}\cos t\cos w)-C_{4}r^{2}\left(  k_{1}^{\prime}\cos t+k_{1}k_{2}\sin
t\right)  \cot w\right)  }{C_{1}r^{3}(1+rk_{1}\cos t\cos w)^{3}}.
\end{array}
\right.
\]
When we equate the functions $\mathfrak{m}(s,t,w)$ found above to each other,
we arrive at the following equations:%
\begin{align}
&  (C_{3}-C_{4}\sin t\cot w)\left(  k_{1}^{\prime}\cot t+k_{1}k_{2}\right)
=0,\label{sari1T11}\\
&  {\small k}_{1}{\small (C}_{1}\sec{\small t}\sec{\small w+k}_{2}%
r{\small (C}_{2}\tan{\small t-C}_{4}\sin{\small t}\cot{\small w))}%
+{\small r(C}_{1}{\small k}_{1}^{2}-{\small k}_{1}^{\prime}{\small (C}_{4}%
\cos{\small t}\cot{\small w-C}_{2}{\small ))=0,}\label{sari2T11}\\
&  k_{1}(rk_{2}(C_{3}-C_{2}\tan t)-C_{1}\sec t\sec w)-C_{1}rk_{1}^{2}%
+rk_{1}^{\prime}(C_{3}\cot t-C_{2})=0. \label{sari3T11}%
\end{align}
In the equation (\ref{sari1T11}), it holds that $k_{1}^{\prime}\cot
t+k_{1}k_{2}\neq0$. This is because, when $k_{1}^{\prime}\cot t+k_{1}k_{2}=0,$
the function $\mathfrak{n}(s,t,w)$ in the first equation of
(\ref{denklemlerT11}) becomes zero. This, in turn, contradicts the definition
of the function $\mathfrak{n}(s,t,w)$ in our classification as $\mathcal{L}%
_{1}N=\mathfrak{m}N+\mathfrak{n}C$. So, from the equation (\ref{sari1T11}) and
$k_{1}^{\prime}\cot t+k_{1}k_{2}\neq0,$ we have $C_{3}=C_{4}=0.$ When
$C_{3}=C_{4}=0$, substituting this into the equation (\ref{sari3T11}) yields%
\[
\left(  C_{1}rk_{1}^{2}+C_{2}rk_{1}^{\prime}\right)  \cos t+C_{1}k_{1}\sec
w+C_{2}rk_{1}k_{2}\sin t=0.
\]
Thus, we have%
\[
C_{1}k_{1}^{2}+C_{2}k_{1}^{\prime}=C_{1}k_{1}=C_{2}k_{1}k_{2}=0
\]
and so $C_{1}=C_{2}=0.$ This is a contradiction.

Secondly, let us assume that $C_{1}=0$.

In this case, from the first equation of the set of equations (\ref{turevT11})
it's easy to see that%
\begin{equation}
C_{2}k_{1}=0{.} \label{C2k1=0}%
\end{equation}
If $k_{1}=0$ in (\ref{C2k1=0}), then from the second, third and fourth
equations of (\ref{denklemlerT11}), it is calculated as%

\begin{equation}
\left.
\begin{array}
[c]{l}%
C_{2}r^{3}\mathfrak{n}(s,t,w)=\left(  \mathfrak{m}(s,t,w)r^{3}-2\right)  \cos
t\cos w,\\
\\
C_{3}r^{3}\mathfrak{n}(s,t,w)=\left(  \mathfrak{m}(s,t,w)r^{3}-2\right)  \sin
t\cos w,\\
\\
C_{4}r^{3}\mathfrak{n}(s,t,w)=\left(  \mathfrak{m}(s,t,w)r^{3}-2\right)  \sin
w,
\end{array}
\right\}  \label{nCi}%
\end{equation}
respectively. Since the functions $C_{i}$ depend only on $s$, there is no
solution for functions $\mathfrak{n}(s,t,w)$ in (\ref{nCi}).

Now, let us assume that $C_{2}=0$ in (\ref{C2k1=0}). In this case, from the
second equation of (\ref{denklemlerT11}), it's easy to see that%
\begin{equation}
\mathfrak{m}(s,t,w)=\frac{2\left(  \cos t\cos w+rk_{1}\left(  \cos
(2t)+3rk_{1}\cos^{3}t\cos^{3}w+2\cos^{2}t\cos(2w)\right)  \right)  }%
{r^{3}(1+rk_{1}\cos t\cos w)^{2}\cos t\cos w}{.} \label{c2esit0nT11}%
\end{equation}
Here, when the equation (\ref{c2esit0nT11}) is successively substituted into
the third and fourth equations of (\ref{denklemlerT11}), we obtain%
\[
\left.
\begin{array}
[c]{l}%
\mathfrak{n}(s,t,w)C_{3}=\frac{-2k_{1}}{r^{2}\left(  1+r\cos t\cos
wk_{1}\right)  ^{2}}\tan t,\\
\\
\mathfrak{n}(s,t,w)C_{4}=\frac{-2k_{1}}{r^{2}\left(  1+rk_{1}\cos t\cos
w\right)  ^{2}}\sec t\tan w.
\end{array}
\right.
\]
Here, there is no solution for functions $\mathfrak{n}(s,t,w)$.

Hence, we can state the following theorem:

\begin{theorem}
There are no tubular hypersurfaces (\ref{surfT11}), obtained by pseudo
hyperspheres whose centers lie on a timelike curve in $\mathbb{E}%
^{4}_{1}$, with generalized $\mathcal{L}_{1}$ 1-type
Gauss map.
\end{theorem}

Now, let us classificate the tubular hypersurfaces $\mathcal{T}(s,t,w)$ which
have second kind $\mathcal{L}_{1}$-pointwise 1-type Gauss map, i.e.,
$\mathcal{L}_{1}N=\mathfrak{m}\left(  N+C\right)  .$

From (\ref{normalT11}) and (\ref{L1T11}), we get%
\begin{equation}
\left.
\begin{array}
[c]{l}%
-\frac{2\left(  k_{1}^{\prime}\cos t+k_{1}k_{2}\sin t\right)  \cos
w}{r(1+rk_{1}\cos t\cos w)^{3}}=\mathfrak{m}C_{1},\\
\\
-\frac{2\left(  rk_{1}\left(  3rk_{1}\cos^{3}t\cos^{3}w+2\cos^{2}%
t\cos(2w)+\cos(2t)\right)  +\cos t\cos w\right)  }{r^{3}(1+rk_{1}\cos t\cos
w)^{2}}=\mathfrak{m}\left(  {-\cos t\cos w+}C_{2}\right)  ,\\
\\
-\frac{2(1+3rk_{1}\cos t\cos w)\sin t\cos w}{r^{3}(1+rk_{1}\cos t\cos
w)}=\mathfrak{m}\left(  -{\sin t\cos w+}C_{3}\right)  ,\\
\\
-\frac{2(1+3rk_{1}\cos t\cos w)\sin w}{r^{3}(1+rk_{1}\cos t\cos w)}%
=\mathfrak{m}\left(  -{\sin w+}C_{4}\right)  .
\end{array}
\right\}  \label{denklerT11}%
\end{equation}
Here, from the fourth equation of (\ref{denklerT11}) it's easy to see that%
\begin{equation}
\mathfrak{m}(s,t,w)={-\frac{2\left(  1+3rk_{1}\cos t\cos w\right)  \sin
w}{r^{3}(1+rk_{1}\cos t\cos w)\left(  -{\sin w+}C_{4}\right)  }.}
\label{4.denkmT11}%
\end{equation}

When the equation (\ref{4.denkmT11}) is successively substituted into the
second and third equations of (\ref{denklerT11}), we obtain%
\begin{align*}
&  3C_{2}r^{2}k_{1}^{2}\cos^{2}t\cos^{2}w+4C_{2}rk_{1}\cos t\cos
w+C_{2}-rk_{1}=0,\\
&  C_{4}\sin t\cos w-C_{3}\sin w=0.
\end{align*}
So, we have%
\begin{equation}
k_{1}=C_{2}=C_{3}=C_{4}=0. \label{denk1T11}%
\end{equation}
Now, when the components of the equation (\ref{denk1T11}) is substituted into
the second or third equations of (\ref{denklerT11}), we calculated%
\begin{equation}
\mathfrak{m}(s,t,w)={\frac{2}{r^{3}}.} \label{denkm2/r3T11}%
\end{equation}
Also, from the first equation of (\ref{denklerT11}) and (\ref{denkm2/r3T11}),
we have $C_{1}=0.$

From the calculations made above for classificate the tubular hypersurfaces
$\mathcal{T}(s,t,w)$ which have second kind $\mathcal{L}_{1}$-pointwise 1-type
Gauss map, i.e., $\mathcal{L}_{1}N=\mathfrak{m}\left(  N+C\right)  $, we can
give the following theorem:

\begin{theorem}
There are no tubular hypersurfaces (\ref{surfT11}), obtained by pseudo
hyperspheres whose centers lie on a timelike curve in $\mathbb{E}%
^{4}_{1}$, with second kind $\mathcal{L}_{1}$-pointwise
1-type Gauss map.
\end{theorem}

Moreover, if the function $m$ is constant in Definition \ref{def} (\textbf{ii
}or \textbf{iii}), then we say $M$ has first or second kind $\mathcal{L}_{k}%
$-(global) pointwise 1-type Gauss map. Thus, we can state the following theorem:

\begin{theorem}
The tubular hypersurfaces (\ref{surfT11}), obtained by pseudo hyperspheres
whose centers lie on a timelike curve in $\mathbb{E}%
^{4}_{1}$, have first kind $\mathcal{L}_{1}$-(global)
pointwise 1-type Gauss map, i.e., $\mathcal{L}_{1}N=\mathfrak{m}N$ if and only if $k_{1}=0,$ where $\mathfrak{m}(s,t,w)={\frac
{2}{r^{3}}.}$
\end{theorem}

Finally, in the equation (\ref{L1T11}), since the coefficients of
$F_{1},~F_{2},~F_{3}$ and $F_{4}$ cannot all be zero, we can give the
following theorem:

\begin{theorem}
The tubular hypersurfaces (\ref{surfT11}), obtained by pseudo hyperspheres
whose centers lie on a timelike curve in $\mathbb{E}%
^{4}_{1}$, cannot have $\mathcal{L}_{1}$-harmonic Gauss map.
\end{theorem}

\subsection{Some Classifications for Tubular Hypersurfaces Generated by Timelike Curves with $\mathcal{L}_{2}$ Operator in $\mathbb{E}%
^{4}_{1}$}

\

Firstly, it is calculated from (\ref{gradf}), (\ref{gT11}) and (\ref{akT11})
as%
\begin{align}
\nabla a_{3}=  &  \frac{\left(  k_{1}^{\prime}\cos t+k_{1}k_{2}\sin t\right)
\cos w}{r^{2}(1+rk_{1}\cos t\cos w)^{3}}{}F_{1}+\frac{k_{1}\left(  2\cos
^{2}t\cos(2w)+\cos(2t)-3\right)  }{4r^{3}(1+rk_{1}\cos t\cos w)^{2}}{}%
F_{2}\nonumber\\
&  +\frac{k_{1}\sin t\cos t\cos^{2}w}{r^{3}(1+rk_{1}\cos t\cos w)^{2}}{}%
F_{3}+\frac{k_{1}\cos t\sin w\cos w}{r^{3}(1+rk_{1}\cos t\cos w)^{2}}{}F_{4}.
\label{a3nablaT11}%
\end{align}

So, from (\ref{Lk}), (\ref{Hk}), (\ref{normalT11}), (\ref{epsilonT11}),
(\ref{akT11}) and (\ref{a3nablaT11}), we have%
\begin{align}
\mathcal{L}_{2}N  &  =\frac{\left(  k_{1}^{\prime}\cos t+k_{1}k_{2}\sin
t\right)  \cos w}{r^{2}(1+rk_{1}\cos t\cos w)^{3}}{}F_{1}\nonumber\\
&  +\frac{k_{1}\left(
\begin{array}
[c]{l}%
\frac{rk_{1}}{2}\left(
\begin{array}
[c]{l}%
24rk_{1}\cos^{4}t\cos^{4}w+12\cos^{3}t\cos(3w)\\
+19\cos t\cos w+9\cos(3t)\cos w
\end{array}
\right) \\
+6\cos^{2}t\cos(2w)+3\cos(2t)-1
\end{array}
\right)  }{4r^{3}(1+rk_{1}\cos t\cos w)^{3}}{}F_{2}\nonumber\\
&  +\frac{3k_{1}\sin t\cos t\cos^{2}w}{r^{3}(1+rk_{1}\cos t\cos w)}{}%
F_{3}+\frac{3k_{1}\cos t\sin w\cos w}{r^{3}(1+rk_{1}\cos t\cos w)}{}F_{4}.
\label{L2T11}%
\end{align}

Now, let us give some classifications for the tubular hypersurfaces
(\ref{surfT11}) which have generalized $\mathcal{L}_{2}$ 1-type Gauss
map,\ first kind $\mathcal{L}_{2}$-pointwise 1-type Gauss map and second kind
$\mathcal{L}_{2}$-pointwise 1-type Gauss map and $\mathcal{L}_{2}$-harmonic Gauss map.

Now, let us classificate the tubular hypersurfaces $\mathcal{T}(s,t,w)$ which
have generalized $\mathcal{L}_{2}$ 1-type Gauss map. From (\ref{normalT11})
and (\ref{L2T11}), we get%
\begin{equation}
\left.
\begin{array}
[c]{l}%
\frac{\left(  k_{1}^{\prime}\cos t+k_{1}k_{2}\sin t\right)  \cos w}%
{r^{2}(1+rk_{1}\cos t\cos w)^{3}}=\mathfrak{n}C_{1},\\
\\
\frac{k_{1}\left(
\begin{array}
[c]{l}%
\frac{rk_{1}}{2}\left(
\begin{array}
[c]{l}%
24rk_{1}\cos^{4}t\cos^{4}w+12\cos^{3}t\cos(3w)\\
+19\cos t\cos w+9\cos(3t)\cos w
\end{array}
\right) \\
+6\cos^{2}t\cos(2w)+3\cos(2t)-1
\end{array}
\right)  }{4r^{3}(1+rk_{1}\cos t\cos w)^{3}}=\mathfrak{m}\left(  {-\cos t\cos
w}\right)  +\mathfrak{n}C_{2},\\
\\
\frac{3k_{1}\sin t\cos t\cos^{2}w}{r^{3}(1+rk_{1}\cos t\cos w)}=\mathfrak{m}%
\left(  -{\sin t\cos w}\right)  +\mathfrak{n}C_{3},\\
\\
\frac{3k_{1}\cos t\sin w\cos w}{r^{3}(1+rk_{1}\cos t\cos w)}=\mathfrak{m}%
\left(  -{\sin w}\right)  +\mathfrak{n}C_{4}.
\end{array}
\right\}  \label{L2.1.denklemler}%
\end{equation}
Firstly, let us assume that $C_{1}\neq0$.

In this case, from the first equation of (\ref{L2.1.denklemler}) it's easy to
see that%
\begin{equation}
\mathfrak{n}(s,t,w)=\frac{\left(  k_{1}^{\prime}\cos t+k_{1}k_{2}\sin
t\right)  \cos w}{r^{2}(1+rk_{1}\cos t\cos w)^{3}C_{1}}{.}
\label{L2c1esitdegil0nT11}%
\end{equation}
Here, when the equation (\ref{L2c1esitdegil0nT11}) is successively substituted
into the second, third and fourth equations of (\ref{L2.1.denklemler}), we
obtain%
\[
\left.
\begin{array}
[c]{l}%
\mathfrak{m}(s,t,w)=\frac{\left(
\begin{array}
[c]{l}%
2k_{1}(-3C_{1}\cos t\cos w+C_{1}\sec t\sec w+C_{2}rk_{2}\tan t)+2C_{2}%
rk_{1}^{\prime}\\
-6C_{1}r^{2}k_{1}^{3}\cos^{3}t\cos^{3}w-C_{1}rk_{1}^{2}\left(  6\cos^{2}%
t\cos(2w)+3\cos(2t)+1\right)
\end{array}
\right)  }{2C_{1}r^{3}(1+rk_{1}\cos t\cos w)^{3}},\\
\\
\mathfrak{m}(s,t,w)=\frac{C_{3}r\left(  k_{1}^{\prime}\cot t+k_{1}%
k_{2}\right)  -3C_{1}k_{1}(1+rk_{1}\cos t\cos w)^{2}\cos t\cos w}{C_{1}%
r^{3}(1+rk_{1}\cos t\cos w)^{3}},\\
\\
\mathfrak{m}(s,t,w)=\frac{C_{4}r\left(  k_{1}^{\prime}\cos t+k_{1}k_{2}\sin
t\right)  \cot w-3C_{1}k_{1}(1+rk_{1}\cos t\cos w)^{2}\cos t\cos w}{C_{1}%
r^{3}(1+rk_{1}\cos t\cos w)^{3}}.
\end{array}
\right.
\]
When we equate the functions $\mathfrak{m}(s,t,w)$ found above to each other,
we arrive at the following equations:%
\begin{align}
&  (C_{4}\sin t\cos w-C_{3}\sin w)\left(  k_{1}^{\prime}\cos t+k_{1}k_{2}\sin
t\right)  =0,\label{L2.sari1T11}\\
&  k_{1}(C_{1}\sec t\sec w+rk_{2}(C_{2}\tan t-C_{4}\sin t\cot w))+C_{1}%
rk_{1}^{2}+rk_{1}^{\prime}(C_{2}-C_{4}\cos t\cot w){\small =0,}%
\label{L2.sari2T11}\\
&  k_{1}(C_{1}\sec t\sec w+rk_{2}(C_{2}\tan t-C_{3}))+C_{1}rk_{1}^{2}%
+rk_{1}^{\prime}(C_{2}-C_{3}\cot t)=0. \label{L2.sari3T11}%
\end{align}
In the equation (\ref{L2.sari1T11}), it holds that $k_{1}^{\prime}\cos
t+k_{1}k_{2}\sin t\neq0$. This is because when $k_{1}^{\prime}\cos
t+k_{1}k_{2}\sin t=0,$ the function $\mathfrak{n}(s,t,w)$ in the first
equation of (\ref{L2.1.denklemler}) becomes zero. This, in turn, contradicts
the definition of the function $\mathfrak{n}(s,t,w)$ in our classification as
$\mathcal{L}_{2}N=\mathfrak{m}N+\mathfrak{n}C$. So, from the equation
(\ref{L2.sari1T11}) and $k_{1}^{\prime}\cos t+k_{1}k_{2}\sin t\neq0,$ we have
$C_{3}=C_{4}=0.$ When $C_{3}=C_{4}=0$, substituting this into the equation
(\ref{L2.sari3T11}) yields%
\[
r\left(  C_{1}k_{1}^{2}+C_{2}k_{1}^{\prime}\right)  \cos t+C_{2}rk_{1}%
k_{2}\sin t+C_{1}k_{1}\sec w=0.
\]
Thus, we have%
\[
C_{1}k_{1}^{2}+C_{2}k_{1}^{\prime}=C_{1}k_{1}=C_{2}k_{1}k_{2}=0
\]
and so $C_{1}=C_{2}=0.$ This is a contradiction.

Secondly, let us assume that $C_{1}=0$.

In this case, from the first equation of the set of equations (\ref{turevT11})
it's easy to see that%
\begin{equation}
C_{2}k_{1}=0{.} \label{L2.C2k1}%
\end{equation}
If $k_{1}=0$ in (\ref{L2.C2k1}), then from the second, third and fourth
equations of (\ref{L2.1.denklemler}), it is calculated as%

\begin{equation}
\left.
\begin{array}
[c]{l}%
\mathfrak{m}(s,t,w)\cos t\cos w=\mathfrak{n}(s,t,w)C_{2},\\
\mathfrak{m}(s,t,w)\sin t\cos w=\mathfrak{n}(s,t,w)C_{3},\\
\mathfrak{m}(s,t,w)\sin w=\mathfrak{n}(s,t,w)C_{4},
\end{array}
\right\}  \label{L2.m+n=0}%
\end{equation}
respectively. Since the functions $C_{i}$ depend only on $s$, there is no
solution for functions $\mathfrak{m}(s,t,w)$ and $\mathfrak{n}(s,t,w)$ in
(\ref{L2.m+n=0}).

Now, let us assume that $C_{2}=0$ in (\ref{L2.C2k1}). In this case, from the
second equation of (\ref{L2.1.denklemler}), it's easy to see that%
\begin{equation}
\mathfrak{m}(s,t,w)=-\frac{k_{1}\left(
\begin{array}
[c]{l}%
rk_{1}\left(
\begin{array}
[c]{l}%
24rk_{1}\cos^{4}t\cos^{4}w+12\cos^{3}t\cos(3w)\\
+19\cos t\cos w+9\cos(3t)\cos w
\end{array}
\right) \\
+12\cos^{2}t\cos(2w)+6\cos(2t)-2
\end{array}
\right)  }{8r^{3}(1+rk_{1}\cos t\cos w)^{3}\cos t\cos w}{.}
\label{L2.c2esit0nT11}%
\end{equation}
Here, when the equation (\ref{L2.c2esit0nT11}) is successively substituted
into the third and fourth equations of (\ref{L2.1.denklemler}), we obtain%
\[
\left.
\begin{array}
[c]{l}%
\mathfrak{n}(s,t,w)C_{3}=\frac{k_{1}(rk_{1}\sin t\cos w+\tan t)}%
{r^{3}(1+rk_{1}\cos t\cos w)^{3}},\\
\\
\mathfrak{n}(s,t,w)C_{4}=\frac{k_{1}(rk_{1}\sin w+\sec t\tan w)}%
{r^{3}(1+rk_{1}\cos t\cos w)^{3}}.
\end{array}
\right.
\]
Here, there is no solution for functions $\mathfrak{n}(s,t,w)$.

Therefore, we can give the following theorem:

\begin{theorem}
There are no tubular hypersurfaces (\ref{surfT11}), obtained by pseudo
hyperspheres whose centers lie on a timelike curve in $\mathbb{E}%
^{4}_{1}$, with generalized $\mathcal{L}_{2}$ 1-type
Gauss map.
\end{theorem}

Now, let us classificate the tubular hypersurfaces $\mathcal{T}(s,t,w)$ which
have second kind $\mathcal{L}_{2}$-pointwise 1-type Gauss map, i.e.,
$\mathcal{L}_{2}N=\mathfrak{m}\left(  N+C\right)  .$

From (\ref{normalT11}) and (\ref{L2T11}), we get%
\begin{equation}
\left.
\begin{array}
[c]{l}%
\frac{\left(  k_{1}^{\prime}\cos t+k_{1}k_{2}\sin t\right)  \cos w}%
{r^{2}(1+rk_{1}\cos t\cos w)^{3}}=\mathfrak{m}C_{1},\\
\\
\frac{k_{1}\left(
\begin{array}
[c]{l}%
\frac{rk_{1}}{2}\left(
\begin{array}
[c]{l}%
24rk_{1}\cos^{4}t\cos^{4}w+12\cos^{3}t\cos(3w)\\
+19\cos t\cos w+9\cos(3t)\cos w
\end{array}
\right) \\
+6\cos^{2}t\cos(2w)+3\cos(2t)-1
\end{array}
\right)  }{4r^{3}(1+rk_{1}\cos t\cos w)^{3}}=\mathfrak{m}\left(  {-\cos t\cos
w+}C_{2}\right)  ,\\
\\
\frac{3k_{1}\sin t\cos t\cos^{2}w}{r^{3}(1+rk_{1}\cos t\cos w)}=\mathfrak{m}%
\left(  -{\sin t\cos w+}C_{3}\right)  ,\\
\\
\frac{3k_{1}\cos t\sin w\cos w}{r^{3}(1+rk_{1}\cos t\cos w)}=\mathfrak{m}%
\left(  -{\sin w+}C_{4}\right)  .
\end{array}
\right\}  \label{L2.denkT11ler}%
\end{equation}
Here, from the last equation of (\ref{L2.denkT11ler}) it's easy to see that%
\begin{equation}
\mathfrak{m}(s,t,w)=\frac{3k_{1}\cos t\sin w\cos w}{r^{3}(1+rk_{1}\cos t\cos
w)\left(  -{\sin w+}C_{4}\right)  }{.} \label{L2.4.denkmT11}%
\end{equation}
Here, when the equation (\ref{L2.4.denkmT11}) is substituted into the second
equation of (\ref{L2.denkT11ler}), we obtain%
\[
-1+3C_{2}\cos t\cos w+3C_{2}rk_{1}\cos^{2}t\cos^{2}w=0.
\]
Since the last equality is never zero, we can give the following theorem:

\begin{theorem}
There are no tubular hypersurface (\ref{surfT11}), obtained by pseudo
hyperspheres whose centers lie on a timelike curve in $\mathbb{E}%
^{4}_{1}$, with second kind $\mathcal{L}_{2}$-pointwise
1-type Gauss map.
\end{theorem}

Now, let us classificate the tubular hypersurfaces $\mathcal{T}(s,t,w)$ which
have first kind $\mathcal{L}_{2}$-pointwise 1-type Gauss map, i.e.,
$\mathcal{L}_{2}N=\mathfrak{m}N.$

From (\ref{normalT11}) and (\ref{L2T11}), we get%
\begin{equation}
\left.
\begin{array}
[c]{l}%
\frac{\cos w\left(  \cos tk_{1}^{\prime}+k_{1}k_{2}\sin t\right)  }%
{r^{2}(1+rk_{1}\cos t\cos w)^{3}}=0,\\
\\
\frac{k_{1}\left(
\begin{array}
[c]{l}%
\frac{1}{2}rk_{1}\left(
\begin{array}
[c]{l}%
24rk_{1}\cos^{4}t\cos^{4}w+12\cos^{3}t\cos(3w)\\
+19\cos t\cos w+9\cos(3t)\cos w
\end{array}
\right) \\
+6\cos^{2}t\cos(2w)+3\cos(2t)-1
\end{array}
\right)  }{4r^{3}(1+rk_{1}\cos t\cos w)^{3}}=\mathfrak{m}\left(  {-\cos t\cos
w}\right)  ,\\
\\
\frac{3k_{1}\sin t\cos t\cos^{2}w}{r^{3}(1+rk_{1}\cos t\cos w)}=\mathfrak{m}%
\left(  -{\sin t\cos w}\right)  ,\\
\\
\frac{3k_{1}\cos t\sin w\cos w}{r^{3}(1+rk_{1}\cos t\cos w)}=\mathfrak{m}%
\left(  -{\sin w}\right)  .
\end{array}
\right\}  \label{L2.denkmNT11ler}%
\end{equation}
Here, from the last equation of (\ref{L2.denkmNT11ler}) it's easy to see that%
\begin{equation}
\mathfrak{m}(s,t,w)=\frac{-3k_{1}\cos t\cos w}{r^{3}(1+rk_{1}\cos t\cos w)}{.}
\label{L2.4.denkmNT11}%
\end{equation}
Here, when the equation (\ref{L2.4.denkmNT11}) is substituted into the second
equation of (\ref{L2.denkT11ler}), we obtain%
\[
k_{1}(rk_{1}+\sec t\sec w)=0.
\]
Since the last equality is never zero, we can give the following theorem:

\begin{theorem}
There are no tubular hypersurfaces (\ref{surfT11}), obtained by pseudo
hyperspheres whose centers lie on a timelike curve in $\mathbb{E}%
^{4}_{1}$, with first kind $\mathcal{L}_{2}$-pointwise
1-type Gauss map.
\end{theorem}

Finally, since the coefficients $F_{1},~F_{2},~F_{3}$ and $F_{4}$ in equation
(\ref{L2T11}) are all zero only for $k_{1}=0$, we can give the following theorem:

\begin{theorem}
The tubular hypersurfaces (\ref{surfT11}), obtained by pseudo hyperspheres
whose centers lie on a timelike curve in $\mathbb{E}%
^{4}_{1}$, have $\mathcal{L}_{2}$-harmonic Gauss map if and only if
$k_{1}=0$.
\end{theorem}

\section{{Some Classifications for Tubular Hypersurfaces Generated by
Spacelike Curves with }$\mathcal{L}_{k}$ Operators in $\mathbb{E}%
^{4}_{1}$}

In this section, we give the general formulas for $\mathcal{L}_{1}$
(Cheng-Yau) and $\mathcal{L}_{2}$ operators of the Gauss maps of the six types
of tubular hypersurfaces $\mathcal{T}^{\{j,\lambda\}}(s,t,w)$ that are formed
as the envelope of a family of pseudo hyperspheres{ or pseudo hyperbolic
hyperspheres} whose centers lie on spacelike curves $\beta(s)$ with non-null
Frenet vectors in $\mathbb{E}%
^{4}_{1}$ and give some classifications for these
hypersurfaces which have generalized $\mathcal{L}_{k}$ 1-type Gauss map, first
kind $\mathcal{L}_{k}$-pointwise 1-type Gauss map and second kind
$\mathcal{L}_{k}$-pointwise 1-type Gauss map and $\mathcal{L}_{k}%
$-harmonic Gauss map, $k\in\{1,2\}$.

The tubular hypersurfaces $\mathcal{T}^{\{j,\lambda\}}(s,t,w)$ that are formed
as the envelope of a family of pseudo hyperspheres or pseudo hyperbolic
hyperspheres whose centers lie on a spacelike curve with non-null Frenet
vectors $F_{i}$ in $\mathbb{E}%
^{4}_{1}$ can be parametrized by%
\begin{equation}
\left.
\begin{array}
[c]{l}%
\mathcal{T}^{\{2,1\}}(s,t,w)=\beta(s)+r(\cosh t\sinh wF_{2}(s)+\cosh
wF_{3}(s)+\sinh t\sinh wF_{4}(s)),\\
\mathcal{T}^{\{2,-1\}}(s,t,w)=\beta(s)+r(\cosh t\cosh wF_{2}(s)+\sinh
wF_{3}(s)+\sinh t\cosh wF_{4}(s)),\\
\mathcal{T}^{\{3,1\}}(s,t,w)=\beta(s)+r(\sinh t\sinh wF_{2}(s)+\cosh t\sinh
wF_{3}(s)+\cosh wF_{4}(s)),\\
\mathcal{T}^{\{3,-1\}}(s,t,w)=\beta(s)+r(\sinh t\cosh wF_{2}(s)+\cosh t\cosh
wF_{3}(s)+\sinh wF_{4}(s)),\\
\mathcal{T}^{\{4,1\}}(s,t,w)=\beta(s)+r(\cosh wF_{2}(s)+\sinh t\sinh
wF_{3}(s)+\cosh t\sinh wF_{4}(s)),\\
\mathcal{T}^{\{4,-1\}}(s,t,w)=\beta(s)+r(\sinh wF_{2}(s)+\sinh t\cosh
wF_{3}(s)+\cosh t\cosh wF_{4}(s)),
\end{array}
\right\}  \label{surfTjlamda}%
\end{equation}
respectively. Here, we suppose for $\mathcal{T}^{\{j;\lambda\}}(s,t,w)$ that

i) $\left\langle F_{j},F_{j}\right\rangle =-1=\varepsilon_{j}$\ and\ for
$i\neq j,$ $\left\langle F_{i},F_{i}\right\rangle =1=\varepsilon_{i},$
$i,j\in\{1,2,3,4\},$

ii) {if the tubular hypersurface is foliated by pseudo hyperspheres or pseudo
hyperbolic hyperspheres, then }$\lambda=1$\textit{ or }$\lambda=-1,${
respectively }(for more details, one can see \cite{BKMSbizim}).

Now, let us write the following lemma which states the general parametric
expressions of 6 different types of tubular hypersurfaces given by
(\ref{surfTjlamda}) and obtained by pseudo hyperspheres and pseudo hyperbolic
hyperspheres whose centers lie on a spacelike curve with non-null Frenet
vector fields in $\mathbb{E}%
^{4}_{1}$.

\begin{lemma}
\label{lemmatoplusurfTjlamda}The general expression of the tubular
hypersurfaces $\mathcal{T}^{\{j,\lambda\}}(s,t,w)$ that are formed as the
envelope of a family of pseudo hyperspheres or pseudo hyperbolic hyperspheres
whose centers lie on a spacelike curve $\beta(s)$ with non-null Frenet vectors
$F_{i}(s)$ in $\mathbb{E}%
^{4}_{1}$ can be given by%
\begin{equation}
\mathcal{T}^{\{j,\lambda\}}(s,t,w)=\beta(s)+r\left(
%TCIMACRO{\tsum _{i=2}^{4}}%
%BeginExpansion
{\textstyle\sum_{i=2}^{4}}
%EndExpansion
\mu_{i}^{\lambda}(s,t,w)F_{i}(s)\right)  , \label{toplusurfTjlamda}%
\end{equation}
where$\ $%
\[
\mu_{5}^{\lambda}(s,t,w)=\mu_{2}^{\lambda}(s,t,w),\text{ }\mu_{6}^{\lambda
}(s,t,w)=\mu_{3}^{\lambda}(s,t,w)
\]
and for $j=2,3,4$
\[
\left.
\begin{array}
[c]{l}%
\mu_{j}^{\lambda}(s,t,w)=\left(  \sinh w\right)  ^{\frac{1+\lambda}{2}}\left(
\cosh w\right)  ^{\frac{1-\lambda}{2}}\cosh t,\\
\\
\mu_{j+1}^{\lambda}(s,t,w)=\left(  \sinh w\right)  ^{\frac{1-\lambda}{2}%
}\left(  \cosh w\right)  ^{\frac{1+\lambda}{2}},\\
\\
\mu_{j+2}^{\lambda}(s,t,w)=\left(  \sinh w\right)  ^{\frac{1+\lambda}{2}%
}\left(  \cosh w\right)  ^{\frac{1-\lambda}{2}}\sinh t.
\end{array}
\right.  \text{ }%
\]
Here, {if the canal hypersurface is foliated by pseudo hyperspheres or pseudo
hyperbolic hyperspheres, then }$\lambda=1$\textit{ or }$\lambda=-1,${
respectively.}
\end{lemma}

Here, we can give the general parametric expressions of the unit normal vector
fields, the coefficients of the first fundamental forms and the principal
curvatures of the tubular hypersurfaces $\mathcal{T}^{\{j,\lambda\}}$
{parametrized by (\ref{toplusurfTjlamda}).}

The unit normal vector fields $N^{\{j,\lambda\}}$ $(j=2,3,4)$ of
{(\ref{toplusurfTjlamda}) are}%
\begin{equation}
N^{\{j,\lambda\}}=-(-1)^{(4-j)!}\lambda^{j}%
%TCIMACRO{\tsum _{i=2}^{4}}%
%BeginExpansion
{\textstyle\sum_{i=2}^{4}}
%EndExpansion
\mu_{i}^{\lambda}F_{i}. \label{toplu Njlamda}%
\end{equation}
The coefficients of the first fundamental forms $g_{ik}^{\{j,\lambda\}}$
$(j=2,3,4)$ of {(\ref{toplusurfTjlamda}) are}%
\begin{equation}
\left.
\begin{array}
[c]{l}%
g_{11}^{\{j,\lambda\}}=1+r^{2}(k_{2})^{2}\left(  -(-1)^{(4-j)!}\left(  \mu
_{3}^{\lambda}\right)  ^{2}+(-1)^{j}\left(  \mu_{2}^{\lambda}\right)
^{2}\right) \\
\text{ \ \ \ \ \ \ \ \ \ }+r^{2}(k_{3})^{2}\left(  (-1)^{(5-j)!}\left(
\mu_{3}^{\lambda}\right)  ^{2}+(-1)^{j}\left(  \mu_{4}^{\lambda}\right)
^{2}\right) \\
\text{ \ \ \ \ \ \ \ \ \ }+2(-1)^{(4-j)!}rk_{1}\mu_{2}^{\lambda}+r^{2}%
(k_{1})^{2}\left(  \mu_{2}^{\lambda}\right)  ^{2}-2(-1)^{(5-j)!}r^{2}%
k_{2}k_{3}\mu_{2}^{\lambda}\mu_{4}^{\lambda},\\
\\
g_{12}^{\{j,\lambda\}}=g_{21}^{\{j,\lambda\}}=r^{2}(\mu_{j+1}^{\lambda}%
)_{w}\left(  (-1)^{j}k_{3}\left(  \mu_{2}^{\lambda}\right)  _{w}%
-(-1)^{(4-j)!}k_{2}\left(  \mu_{4}^{\lambda}\right)  _{w}\right), \\
\\
g_{22}^{\{j,\lambda\}}=r^{2}\left(  (\mu_{j+1}^{\lambda})_{w}\right)  ^{2},\\
\\
g_{13}^{\{2,\lambda\}}=g_{31}^{\{2,\lambda\}}=\lambda r^{2}\left(  -k_{2}\cosh
t+k_{3}\sinh t\right)  ,\\
g_{13}^{\{3,\lambda\}}=g_{31}^{\{3,\lambda\}}=-\lambda r^{2}k_{3}\cosh t,\\
g_{13}^{\{4,\lambda\}}=g_{31}^{\{4,\lambda\}}=\lambda r^{2}k_{2}\sinh t,\\
\\
g_{23}^{\{j,\lambda\}}=g_{32}^{\{j,\lambda\}}=0,\\
\\
g_{33}^{\{j,\lambda\}}=-\lambda r^{2}.%
\end{array}
\right\}  \label{toplugij}%
\end{equation}
The principal curvatures $\kappa_{i}^{\{j,\lambda\}}$ $(j=2,3,4)$ of
{(\ref{toplusurfTjlamda}) are}%
\begin{equation}
\left.
\begin{array}
[c]{l}%
\kappa_{1}^{\{j,\lambda\}}=\kappa_{2}^{\{j,\lambda\}}=\frac{(-1)^{(4-j)!}%
\lambda^{j}}{r},\\
\\
\kappa_{3}^{\{j,\lambda\}}=\frac{k_{1}\mu_{2}^{\lambda}}{\lambda^{j}\left(
1+(-1)^{(4-j)!}rk_{1}\mu_{2}^{\lambda}\right)  }.
\end{array}
\right\}  \label{topluki}%
\end{equation}
From Lemma \ref{lemmatoplusurfTjlamda}, (\ref{toplu Njlamda}), (\ref{toplugij}%
) and (\ref{topluki}), we get%
\begin{align}
\mathcal{L}_{1}N^{\{j,\lambda\}}  &  =\frac{2\left(  -(-1)^{(5-j)!}k_{1}%
k_{2}\mu_{3}^{\lambda}+k_{1}^{\prime}\mu_{2}^{\lambda}\right)  }{r\left(
(-1)^{(4-j)!}+rk_{1}\mu_{2}^{\lambda}\right)  ^{3}}F_{1}\nonumber\\
&  +\frac{-2\lambda\left(  \mu_{2}^{\lambda}+3r^{2}(\mu_{2}^{\lambda}%
)^{3}\left(  k_{1}\right)  ^{2}+k_{1}\left(  \lambda r+4(-1)^{(4-j)!}r(\mu
_{2}^{\lambda})^{2}\right)  \right)  }{r^{3}\left(  (-1)^{(4-j)!}+rk_{1}%
\mu_{2}^{\lambda}\right)  ^{2}}F_{2}\nonumber\\
&  +\frac{-2\lambda(-1)^{(4-j)!}\mu_{3}^{\lambda}\left(  1+3(-1)^{(4-j)!}%
rk_{1}\mu_{2}^{\lambda}\right)  }{r^{3}\left(  (-1)^{(4-j)!}+rk_{1}\mu
_{2}^{\lambda}\right)  }F_{3}+\frac{-2\lambda\mu_{4}^{\lambda}\left(
(-1)^{(4-j)!}+3rk_{1}\mu_{2}^{\lambda}\right)  }{r^{3}\left(  (-1)^{(4-j)!}%
+rk_{1}\mu_{2}^{\lambda}\right)  }F_{4}. \label{L1jlamda}%
\end{align}

Let $\mathcal{T}^{\left\{  j,\lambda\right\}  }(s,t,w)$ have generalized
$L_{1}$ (Cheng-Yau) 1-type Gauss map, i.e., $\mathcal{L}_{1}N^{\{j,\lambda\}}=\mathfrak{m}%
N^{\{j,\lambda\}}+\mathfrak{n}C,$ where $C=C_{1}F_{1}+C_{2}F_{2}+C_{3}F_{3}+C_{4}F_{4}$ is a
constant vector. Here, by taking derivatives of the constant vector $C$ with
respect to $s,$ from (\ref{fr2})-(\ref{fr4}) we obtain for $\mathcal{T}^{\left\{  j,\lambda\right\}  }$ that%
\begin{equation}
\left.
\begin{array}
[c]{l}%
C_{1}^{\prime}+(-1)^{(4-j)!}C_{2}k_{1}=0,\\
C_{2}^{\prime}+C_{1}k_{1}+(-1)^{(5-j)!}C_{3}k_{2}=0,\\
C_{3}^{\prime}+C_{2}k_{2}-(-1)^{(4-j)!}C_{4}k_{3}=0,\\
C_{4}^{\prime}+C_{3}k_{3}=0.
\end{array}
\right\}  \label{Cjlamda}%
\end{equation}
Also, by taking derivatives the constant vector $C$ with respect to$~t,~w$
separately, one can see that the functions $C_{i}$ depend\ only on $s.$

So, with similar procedure in Subsection 3.1, we can give the following theorems:

\begin{theorem}
There are no tubular hypersurfaces (\ref{toplusurfTjlamda}), obtained by pseudo
hyperspheres and pseudo hyperbolic hyperspheres whose centers lie on a
spacelike curve with non-null Frenet vector fields {$F_{i}$ in }$\mathbb{E}%
^{4}_{1}$,
with generalized $\mathcal{L}_{1}$ 1-type Gauss map in $\mathbb{E}%
^{4}_{1}$.
\end{theorem}

\begin{theorem}
There are no tubular hypersurfaces (\ref{toplusurfTjlamda}) obtained by pseudo
hyperspheres and pseudo hyperbolic hyperspheres whose centers lie on a
spacelike curve with non-null Frenet vector fields {$F_{i}$ in }$\mathbb{E}%
^{4}_{1}$
with second kind $\mathcal{L}_{1}$-pointwise 1-type Gauss map in $\mathbb{E}%
^{4}_{1}$.
\end{theorem}

\begin{theorem}
The tubular hypersurfaces (\ref{toplusurfTjlamda}) obtained by pseudo
hyperspheres and pseudo hyperbolic hyperspheres whose centers lie on a
spacelike curve with non-null Frenet vector fields {$F_{i}$ in }$\mathbb{E}%
^{4}_{1}$
{have} first kind $\mathcal{L}_{1}$-(global) pointwise 1-type Gauss map, i.e.,
$\mathcal{L}_{1}N^{\{j,\lambda\}}=\mathfrak{m}N^{\{j,\lambda\}}$ in $\mathbb{E}%
^{4}_{1}$ if and only if $k_{1}=0,$ where
$\mathfrak{m}(s,t,w)={\frac{2\lambda^{j+1}(-1)^{(4-j)!}}{r^{3}}.}$
\end{theorem}

\begin{theorem}
The tubular hypersurfaces (\ref{toplusurfTjlamda}) obtained by pseudo
hyperspheres and pseudo hyperbolic hyperspheres whose centers lie on a
spacelike curve with non-null Frenet vector fields {$F_{i}$ in }$\mathbb{E}%
^{4}_{1}$
cannot have $\mathcal{L}_{1}$-harmonic Gauss map.
\end{theorem}

Also, from Lemma \ref{lemmatoplusurfTjlamda}, (\ref{toplu Njlamda}),
(\ref{toplugij}) and (\ref{topluki}), we get%
\begin{align}
\mathcal{L}_{2}N^{\{j,\lambda\}}  &  =\frac{\lambda^{j}\left(  (-1)^{j}\mu
_{3}^{\lambda}k_{1}k_{2}-(-1)^{(4-j)!}\mu_{2}^{\lambda}k_{1}^{\prime}\right)
}{r^{2}\left(  (-1)^{(4-j)!}+rk_{1}\mu_{2}^{\lambda}\right)  ^{3}}%
F_{1}\nonumber\\
&  +\frac{-\lambda^{j+1}k_{1}\left(  2\lambda(-1)^{(4-j)!}-3(-1)^{j}\left(
\mu_{4}^{\lambda}\right)  ^{2}-3(-1)^{(5-j)!}\left(  \mu_{3}^{\lambda}\right)
^{2}-3(-1)^{(4-j)!}rk_{1}\left(  \mu_{2}^{\lambda}\right)  ^{3}\right)
}{r^{3}\left(  (-1)^{(4-j)!}+rk_{1}\mu_{2}^{\lambda}\right)  ^{2}}%
F_{2}\nonumber\\
&  +\frac{\lambda^{j+1}\mu_{3}^{\lambda}\left(  3(-1)^{(5-j)!}rk_{1}\mu
_{2}^{\lambda}\right)  }{r^{4}\left(  (-1)^{(5-j)!}+(-1)^{j}rk_{1}\mu
_{2}^{\lambda}\right)  }F_{3}+\frac{\lambda^{j+1}\mu_{4}^{\lambda}\left(
3(-1)^{(5-j)!}rk_{1}\mu_{2}^{\lambda}\right)  }{r^{4}\left(  (-1)^{(5-j)!}%
+(-1)^{j}rk_{1}\mu_{2}^{\lambda}\right)  }F_{4}. \label{L2jlamda}%
\end{align}

Thus, with similar procedure in Subsection 3.2, we can give the following theorems:

\begin{theorem}
There are no tubular hypersurfaces (\ref{toplusurfTjlamda}) obtained by pseudo
hyperspheres and pseudo hyperbolic hyperspheres whose centers lie on a
spacelike curve with non-null Frenet vector fields {$F_{i}$ in }$\mathbb{E}%
^{4}_{1}$
with generalized $\mathcal{L}_{2}$ 1-type Gauss map in $\mathbb{E}%
^{4}_{1}$.
\end{theorem}

\begin{theorem}
There are no tubular hypersurfaces (\ref{toplusurfTjlamda}) obtained by pseudo
hyperspheres and pseudo hyperbolic hyperspheres whose centers lie on a
spacelike curve with non-null Frenet vector fields {$F_{i}$ in }$\mathbb{E}%
^{4}_{1}$
with second kind $\mathcal{L}_{2}$-pointwise 1-type Gauss map in $\mathbb{E}%
^{4}_{1}$.
\end{theorem}

\begin{theorem}
There are no tubular hypersurfaces (\ref{toplusurfTjlamda}) obtained by pseudo
hyperspheres and pseudo hyperbolic hyperspheres whose centers lie on a
spacelike curve with non-null Frenet vector fields {$F_{i}$ in }$\mathbb{E}%
^{4}_{1}$
with first kind $\mathcal{L}_{2}$-pointwise 1-type Gauss map in $\mathbb{E}%
^{4}_{1}$.
\end{theorem}

\begin{theorem}
The tubular hypersurfaces (\ref{toplusurfTjlamda}) obtained by pseudo
hyperspheres and pseudo hyperbolic hyperspheres whose centers lie on a
spacelike curve with non-null Frenet vector fields {$F_{i}$ in }$\mathbb{E}%
^{4}_{1}$
have $\mathcal{L}_{2}$-harmonic Gauss map if and only if $k_{1}=0$.
\end{theorem}

\end{document}